\documentclass{article}

\usepackage[final]{neurips_2022}
\usepackage[utf8]{inputenc} 
\usepackage[T1]{fontenc}    
\usepackage{hyperref}       
\usepackage{url}            
\usepackage{booktabs}       
\usepackage{amsfonts}       
\usepackage{nicefrac}       
\usepackage{microtype}      
\usepackage{xcolor}         
\usepackage{amsmath}
\usepackage{amssymb}
\usepackage{graphicx}
\usepackage{mathtools}
\usepackage{algorithm}
\usepackage{algorithmic}
\usepackage{bm}
\usepackage{amsthm}
\usepackage{enumitem}

\title{Approximate FW Algorithm with a novel DMO method over Graph-structured Support Set }
\author{Yijian Pan, Hongjiao Qiang}

\begin{document}
\maketitle

\begin{abstract}
  In this project, we reviewed a paper that deals graph-structured convex optimization (GSCO) problem with the approximate Frank-Wolfe (FW) algorithm. We analyzed and re-implemented the original algorithm and introduced some extensions based on that. Then we conducted experiments to compare the results and concluded that our backtracking line-search method effectively reduced the number of iterations, while our new DMO method (Top-g+ optimal visiting) did not make satisfying enough improvements. 
\end{abstract}

\section{Introduction}
The paper we reviewed introduces a kind of approximate Frank-Wolfe (FW) algorithm to solve convex optimization problems over graph-structured support sets where the linear minimization oracle (LMO) cannot be efficiently obtained in general.

\subsection{Problem Statement}
This paper deals with the following graph-structured convex optimization (GSCO) problem

\begin{equation}
\min_{\bm x \in\mathbb{R}^d} f(\bm x),\text{subject to }\bm x\in \mathcal{ D}(C,\mathbb{M}),
\end{equation}

where \(\mathbb{M}:=\{S_1,S_2,...,S_m\}\), a collection of subsets of \([d]\), with \(\cup_iS_i=[d]\), and \( \mathcal{ D}(C,\mathbb{M})  \triangleq \text{conv}\{\bm x:\|\bm x\|_2 \leq C,\text{supp}(\bm x)\in\mathbb{M}\}\) is a convex hull of the graph-structured support set described by \(\mathbb{M}\), which contains a collection of allowed structures of the problem, and $f$ is a convex differentiable function. The support of \(\bm x,\text{i.e, supp}(\bm x) \triangleq\{i:\bm x_i\neq0\}\), encodes the sparsity pattern of $\bm x$, which can be defined by interesting graph structures such as a path, tree, or cluster over an underlying graph. Model \(\mathbb{M}\) describes many interesting scenarios where graph structures serve as a powerful prior.

\begin{figure}[h]
    \centering
    \includegraphics[width=0.2\textwidth]{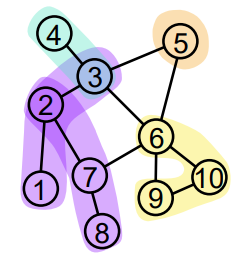}
    \caption{An element of a \(g\)-subgraph model \(\mathbb{M}(\mathbb{G}, s=5, g=4)\)  defined on a 10-node graph where each colored region is a subgraph.}
    \label{figure 1}
\end{figure}

As the authors mentioned, a natural idea to solve the GSCO problem is to use projected gradient descent (PGD) where a projection oracle finds a point in D per per-iteration. In order to obtain approximate convergence guarantees, the existing works of this type generally have the assumption that projection oracles can be solved exactly or with very high approximation guarantees. However, projections satisfying these requirements are usually hard to find for the GSCO problem and multiple projections may be needed at per-iteration. These are the two main issues that the authors look forward to solving with alternative methods.

The authors stated that Frank-Wolfe(FW) methods (Frank et al., 1956), different from PGD-based methods, at each iteration, find a point using the linear minimization oracle (LMO), which for many constraints may enjoy a much cheaper per-iteration cost than the projection oracle (Combettes \& Pokutta, 2021), and often obtain high-quality sparser solutions in early iterations. The fact that the FW methods are less explored for graph-structured optimization problems motivates the authors to tackle the GSCO problem using FW-type methods.

\subsection{Motivation}
The rate of the standard FW method can be improved to \(\mathcal{O}(1/t^2)\) if \(\mathcal{D}\) is a strongly convex set. However, sparsity-including sets are often not strongly convex so the standard FW method is inapplicable to GSCO problems. To achieve the \(\mathcal{O}(1/t^2)\) rate without a strongly convex set, the authors proposed a novel method called approximate FW-type method.

The standard FW method gets a descent direction by using the linear minimization oracle (LMO), which has a cheaper iteration cost than the projection oracle for many constraints. 

\begin{equation}
LMO :  \bm v_t\in \underset{\bm v\in \mathcal{ D}(C,\mathbb{M})}{\mathrm{argmin}}\left<\nabla f(\bm x_t),\bm v \right>
\end{equation}

In general, the LMO returns an extreme point of \(\mathcal{ D}\). Therefore in the case of the \(\mathcal{ D}\) that is defined previously, this problem can be reformulated as a subspace identification problem

\begin{equation}
\underset{\bm v\in \mathcal{ D}(C,\mathbb{M})}{\mathrm{min}}\left<\nabla f(\bm x_t),v \right> = \underset{{S^*}\in \mathbb{M},\|\bm v\|_2 \le \mathcal{ C}}{\mathrm{min}}\left<\nabla f(\bm x_t)_{S^*},\bm v \right>
\end{equation}

where \(S^*\) is an optimal support set minimizing the inner
product in (2). Hence, the minimizer is

\begin{equation}
\bm v_t = - \frac {C \cdot \nabla f(\bm x_t)_{S^*}}{\| \nabla f(\bm x_t)_{S^*} \|_2}, S^* \in \underset{S \in \mathbb{M} }{\mathrm{argmax}} \| \nabla f(\bm x_t)_S \|_2^2 
\end{equation}

The computational complexity of minimizing the inner product (2) increases dramatically due to finding \(S^*\), when \(\mathbb{M}\) is the g-subgraph model or other types of models. To decrease the computational complexity, the authors propose an approximate LMO to find \(\bm v_t\). For our project, in order to extend upon the original authors' work, we planned to modify the process after Top-g+ visiting. In the original algorithm, we arbitrarily added feasible elements to \(\bm v\) so that there would be multiple choices for \(\bm v\). Instead of arbitrarily choosing \(\bm v\) as the proximal gradient, we determined \(\bm v\) by comparing the result of the objective function. The convergence rate can be expected to be improved. 

\section{Related work}

The FW method (Franket et al., 1956) and its variants for convex constrained problems have recently received popularity mainly due to two advantages. First, it is projection free --the LMO is often much cheaper to compute than the projection oracle. Second, in applications with desired structured sparsity, early FW iterations tend to be naturally sparse. Inspired from these advantages, we seek to propose FW-type methods for GSCO problems.

Because of the properties of GSCO problems, using exact FW methods can be impossible. Recent studies mainly focus on the inexact FW methods to solve GSCO problems. The study of inexact FW methods tend to center on two types of LMO errors: gap-additive (Dunn and Harshbarger, 1978; Jaggi, 2013) and gap-multiplicative (Locatello et al., 2017; Pedregosa et al., 2020). Both methods can be hard to solve the LMO without exact support. Gap-additive bound cannot decay properly. Gap-multiplicative estimate could be negative. Instead, rather than approximating the gap, dual maximization oracle (DMO) turns to approximating \(\left \langle \nabla f(\bm x_t),\bm v_t \right \rangle\)(Zhou et al., 2021). An efficient method for DMO to solve GSCO problems is Top-g+ visiting proposed by Zhou. To improve the convergence rate of FW method based on Top-g+ visiting, we propose a new visiting method and combine backtracking linesearch method with standard FW method.

\section{Method}
\subsection{Original Algorithm}
To overcome the computational barrier and improve the convergence rate of FW method, the authors used an inner product operator (IPO) to approximate \(\left \langle \nabla f(\bm x_t),\bm v_t \right \rangle\). Given gradient \(\nabla f(\bm x_t)\), constraint set \(\mathcal{ D}\) and approximation factor \(\delta\in\left( 0,1 \right]\)   the approximated IPO returns \(\bm v_t\) such that 

\begin{equation}
Approximate\; IPO\quad \left<\nabla f(\bm x_t),\bm v_t \right>\le\delta\cdot \min_{\bm s\in\mathcal{ D}}\left<\nabla f(\bm x_t),\bm s \right>
\end{equation}

We denote such set of \(\bm v_t\) as \(\left(\delta,\nabla f(\bm x_t),\mathcal{D}\right)\mathrm{IPO}\).

The approximate IPO can be easier to obtain via DMO. Given the structure support set \(\mathbb{M}\), the DMO finds a set \(S\in\mathbb{M}\) such that 

\begin{equation}
\left \|\nabla f(\bm x_t) _S\right\|_2\ge\delta\cdot\max_{S^\prime \in\mathbb{M}}\left \|\nabla f(\bm x_t) _{S^\prime}\right\|_2
\end{equation}

where approximation factor  \(\delta\in\left( 0,1 \right]\)  . We denote such set \(S\) as the \(\left(\delta,\nabla f(\bm x_t),\mathcal{D}\right)\mathrm{DMO}\).

In practice, the authors used the \(\left(\delta,\nabla f(\bm x_t),\mathcal{D}\right)\mathrm{DMO}\) for the g-subgraph model \(\mathbb{M}\) by visiting Top-g+ neighbors. This DMO algorithm takes in the underlying graph \(\mathbb{G}\), the sparsity \(k\), the number of CCs \(g\), and the gradient vector \(\bm z\) as inputs. It first sorts the entries of \(\bm z\) by magnitudes in order of large to small and picks the \(g\) largest magnitude elements. We initialized the output \(S\) to be the indices of these \(g\) elements in the original vector. Then we numbered these connected components by ID 1 to \(g\). As a next step, we set a graph \(\mathbb{F}\) with edges that were in \(g\) components and then added any feasible additional elements to  \(\mathbb{F}\) and  \(S\) until the number of nonzero elements in \(S\) was equal to \(s\) (the maximum sparsity of \(S\in\mathbb{M}\)). This can be done by just picking any neighbor nodes to the first g “seed” nodes.

Compared with the standard FW method, the approximate FW-type method via DMO has 2 more steps (line 5 and line 6 in Alg. 1) to gain an approximal gradient \(\bm v_t\).

\begin{algorithm}[H]
\caption{FW-type methods for GSCOs}
\begin{algorithmic}[1]
\STATE \textbf{Input:} step size $\{\eta_t\}$, $\delta$, $L$, $C$, and $\mathbb{M}$
\STATE pick any point $\bm x_0$ in $\mathcal{D}$
\FOR{$t = 0,1,\ldots,$}
\STATE $\bm z_t = \begin{cases}\textsc{DMO-FW}:= - \nabla f(\bm x_t) \\
\textsc{DMO-AccFW}:= - \left(\bm x_t - \tfrac{\nabla f(\bm x_t)}{L \eta_t}\right)\end{cases}$
\STATE $S_t = (\delta, -\bm z_t, \mathcal{D})$-DMO
\STATE $\tilde{\bm v}_t =\tfrac{ C \cdot (\bm z_t)_{S_t}}{\left\|( \bm z_t)_{S_t}\right\|_2}$
\STATE Option \textsc{I}:  $\bm x_{t+1} = \bm x_t + \eta_t (\tilde{\bm v}_t - \bm x_t)$
\STATE Option \textsc{II}:  $\bm x_{t+1} = \bm x_t + \eta_t (\tilde{\bm v}_t / \delta - \bm x_t)$ 
\ENDFOR
\STATE \textbf{Return} $(\bm x_{\bar{t}}, f(\bm x_{\bar{t}})), \bar{t} \in \underset{t}{\mathrm{argmin}} f(\bm x_t)$
\end{algorithmic}
\label{algo:approx-fw}
\end{algorithm}

\subsection{Our Extension}
To further improve the rate of approximate FW-type method, we made two changes to the original algorithm. 

First, we added a backtracking line search in determining the optimal step size \(\eta_t\) in every iteration. This could potentially accelerate the convergence of the algorithm, because it ensures that the objective function decreases sufficiently with a properly chosen step size.

The second change we made is optimizing the process of Top-g+ visiting. In the original function, a support set \(S\) is returned whenever it is found. However, multiple latent support sets are available in the constraint set. The support set first found can not guarantee to be the optimal descent direction. To return a better choice of support set \(S\), we compared various feasible support sets \(S\) instead of immediately returning whenever one feasible vector is found. 

We presented the original Top-g+ visiting method in Alg. 3 (in Appendix) while our new method Top-g+ optimal visiting method is in Alg. 4 (in Appendix). There are 2 design differences. First, our new method needs an additional input theta which decides the number of support sets expected to be compared. Second, we modified the process of finding feasible support sets (lines 14-31 in Alg. 3, lines 18-47 in Alg. 4). In our new method, every \((s-g)\) iteration could be viewed as one experiment. After we stayed the \(g\) largest magnitude elements in the dual vector \(z\), we randomly added feasible adjacent nodes to the support set \(S\) until \(\left | {S} \right |=s\) (line 19 in Alg. 4). The process of finding enough feasible nodes to constitute a support set \(S\) needed \((s-g)\) iterations. Once a support set \(S\) was found, we recorded \(S\) and the gradient descent based on the support set (line 25,37 in Alg. 4). After \((s-g) \times \theta\) iterations, we recorded \(\theta\) possible support sets. In the end, we chose the support set \(S\) that had the largest magnitude of gradient descent (line 45,46 in Alg. 4) and returned \(S\). Since the magnitude of gradient descent became larger, the convergence rate could be expected to be improved.

\section{Experiment}

\subsection{Dataset and Optimization Problem}
We finally used the graph from https://github.com/baojian/verse/tree/master/data, and blogcatalog.mat to be specific. There are 10312 nodes, 333983 links, and 39 categories in our dataset. It is composed of 2 sparse matrices, one for the sub-graph information and the other representing the node-to-node connection. Note that the original dataset in the paper is MNIST, which is actually not a graph-like dataset. However, the authors considered it as such for convenience! So we decided to use a real graph-like dataset for the sake of rigor.

The objective function we wanted to optimize is a least-square loss function, i.e.

\begin{equation}
\min_{\bm x \in\mathcal{D} (\mathcal{C},\mathbb{M})} f(\bm x):= \frac{1}{2} \left\|\bm {Ax - y} \right\|_{2}^{2}
\end{equation}

where $\mathcal{D}:= \{\bm x: \operatorname{supp}(\bm x) \in \mathbb{M}, \|\bm x\|_2 \leq C\}$ with $C=1$ and $\mathbb{M} = \{ S\subseteq [d] = [10312]: |S|\leq s=1623\}$. \(\bm A \in \mathbb{R}^{n\times d}\) is a normalized Gaussian sensing matrix where each entry \(a_{ij} \sim N(0,1)\), \(\bm y\) is the observation vector, generated as \(\bm y = \left<\bm A, \bm x^* \right> + \bm e,  \bm e \sim N(0, \sigma^2 \bm I_d)\), and \(\bm x^*\) is a normalized sparse vector with norm 1 that we would like to recover from the observations using the least-square method. We set the number of observations \(n\) to 100 manually.

\subsection{Baseline Implementation}
We first implemented this baseline algorithm with the pseudo-code listed in section 4.2 of the original paper (as well as section 3.1 in our report). The baseline algorithm is a variation of the Frank-Wolfe method, with a DMO method to adapt variables on a graph data set. There are 2 variations with the computation of \(\bm z_t\) in step 4, called DMO-FW and DMO-AccFW respectively, as well as another 2 variations for the update of \(\bm x\) in steps 7 and 8. Besides, we added a stop criterion as $\left |(Obj(x_{t+1})-Obj(x_t)) / Obj(x_t)\right | \le 10^{-6} $. As a default option for the Frank-Wolfe method, the step size $\eta_t$ of the $t$'s iteration here is set as $2/(t+2)$. There are potentially more iteration step size choices, such as exact line search where \(\eta_t = \underset{\eta \ge 0}{\mathrm{argmin}} f(\bm x_t-\eta \nabla f(\bm x_t))\) (which could be expensive to solve), Demyanov-Rubinov step size, where \(\eta_t = min\{ {\frac{\left<-\nabla f(\bm x_t),\tilde{\bm v_t}- \bm x_t\right>}{L \|\tilde{\bm v_t}- \bm x_t\|^2},1}\}\). $L$ is the Lipschitz constant, in the least-square objective function here which can be expressed as $\sigma_1(\bm A^T\bm A)$, which is the largest singular value of \(\bm A^T\bm A\). Because of the time limit and the huge computation cost to get the Lipschitz constant in singular value decomposition, we did not choose the exact line search and Demyanov-Rubinov step size. 

\subsection{Evaluation}
To solve the GSCO problem, a natural idea is to use the projected gradient descent (PGD) method, where a \textit{projection oracle} finds a point in $\mathcal{D}$ at per-iteration. In our Frank-Wolfe method, we can also use the PGD method to find a mask for the variable to satisfy the graph structure. We compared the objective function values for the Frank-Wolfe method via DMO, the random projected gradient descent method, and the best-projected gradient descent by comparing all possible cases (which could be expensive). We didn't implement the DMO-Acc method for it would be impossible to calculate the Lipschitz constant. For all three implementations including the original method and the two improvements, we implemented them based on Option $I$ in Algorithm 1. Here are our results. We can see that the outputs of the Frank-Wolfe method via DMO are very close to the best-projected gradient descent by comparing all possible cases:

\begin{figure}[H]
\centering
\includegraphics[width=.6\textwidth]{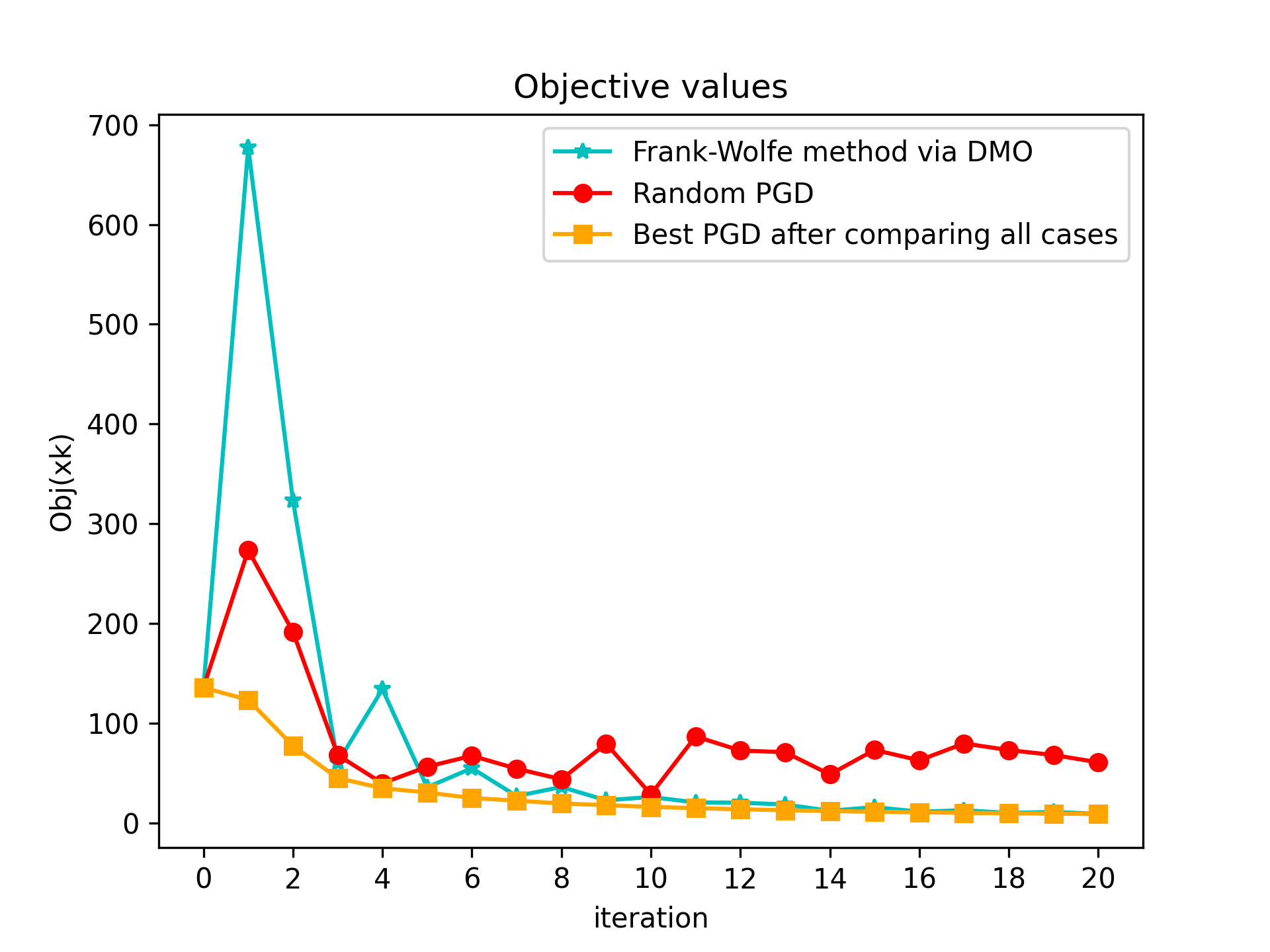}
\caption{Objective function values for the Frank-Wolfe method via DMO, the random projected gradient descent method, and the best-projected gradient descent by comparing all possible cases.}
\end{figure}

\begin{figure}[ht]
\centering
\includegraphics[width=.6\textwidth]{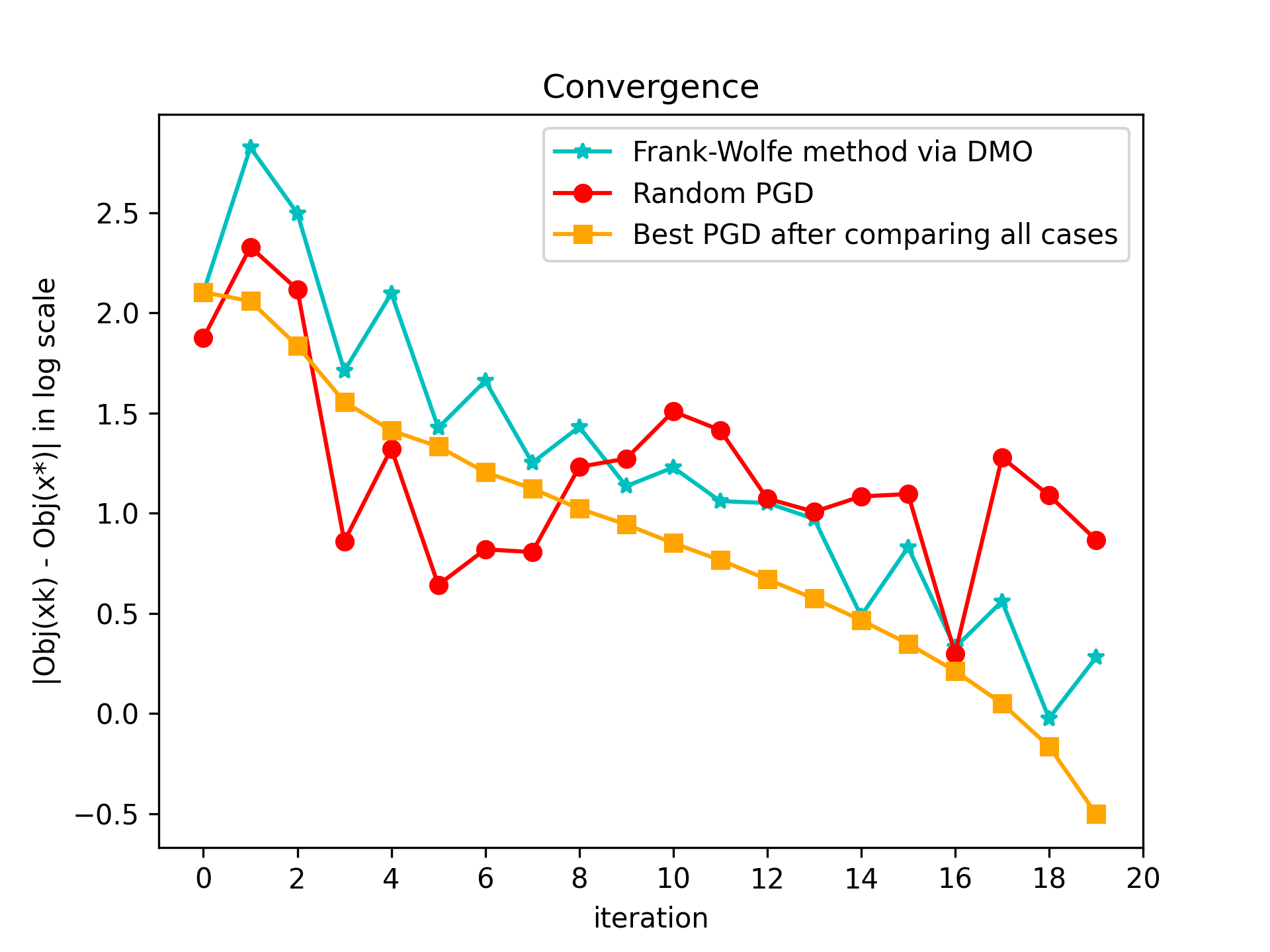}
\caption{Convergence plot for the Frank-Wolfe method via DMO, the random projected gradient descent method, and the best-projected gradient descent by comparing all possible cases.}
\end{figure}

As seen in Fig. 2 and Fig. 3, Best PGD after comparing all the cases has the smallest objection function value in every iteration, and its objective function value steadily decreases as iteration goes. This is anticipated as Best PGD is the most expensive method in computing. Random PGD performs the worst as it has the highest objective function value at the end and the objective function value does not decrease significantly during the iteration process. The objective function value of random PGD shows great fluctuation and even increases significantly at some points. Random PGD projects the gradient on a lower-dimensional subspace randomly, so it may not be as accurate as other algorithms, which can result in slower convergence or lower-quality solutions. Frank-Wolfe method shows similar performance as Best PGD. It achieves a similar final objective function value and generally decreases as iteration goes on. However, we can see some increasing parts of the curve, and where its objective function values are much larger than Best PGD. This is probably due to a too-large step size selection. It seems like it's time to do some improvements.

\subsection{Improvements}

Backtracking line-search (Armijo, L., 1966) is a method used in many iterative algorithms to find an appropriate step size for the next iteration of the optimization process. The idea behind it is to ensure that the objective function decreases sufficiently by iteratively decreasing the step size until a suitable step size is found that satisfies certain conditions.

The traditional backtracking line-search is based on the gradient of current $\bm x_t$. However, we may find the position of $\bm x_t - \bm v_t$ in the updating step $\bm x_{t+1} =\bm  x_t - \eta_t(\bm x_t - \bm v_t)$ very similar to that of $\nabla f(\bm x_t)$ in gradient decent method. So here $\nabla f(\bm x_t)$ were replaced with $\bm x_t - \bm v_t$ in traditional backtracking line-search, and it worked well. The number of iterations required was reduced successfully!

\begin{algorithm}[H]
\caption{Backtracking Line-search for FW-type method}
\begin{algorithmic}[1]
\STATE \textbf{Input:} current location $\bm x_t$, new direction $\bm v_t$, previous step size $\eta_t$, decay parameter $\beta\in (0,1)$
\WHILE {$Obj(\bm x_t - \eta_t(\bm x_t - \bm  v_t)) > Obj(\bm x_t) - 0.5\eta_t\left\|\bm x_t -\bm  v_t\right\|_{2}^{2}$} 
\STATE $\eta_{t+1} \gets \beta\eta_t$
\ENDWHILE
\STATE \textbf{Return} $\eta_t$
\end{algorithmic}
\end{algorithm}

\begin{figure}[H]
\centering
\includegraphics[width=.6\textwidth]{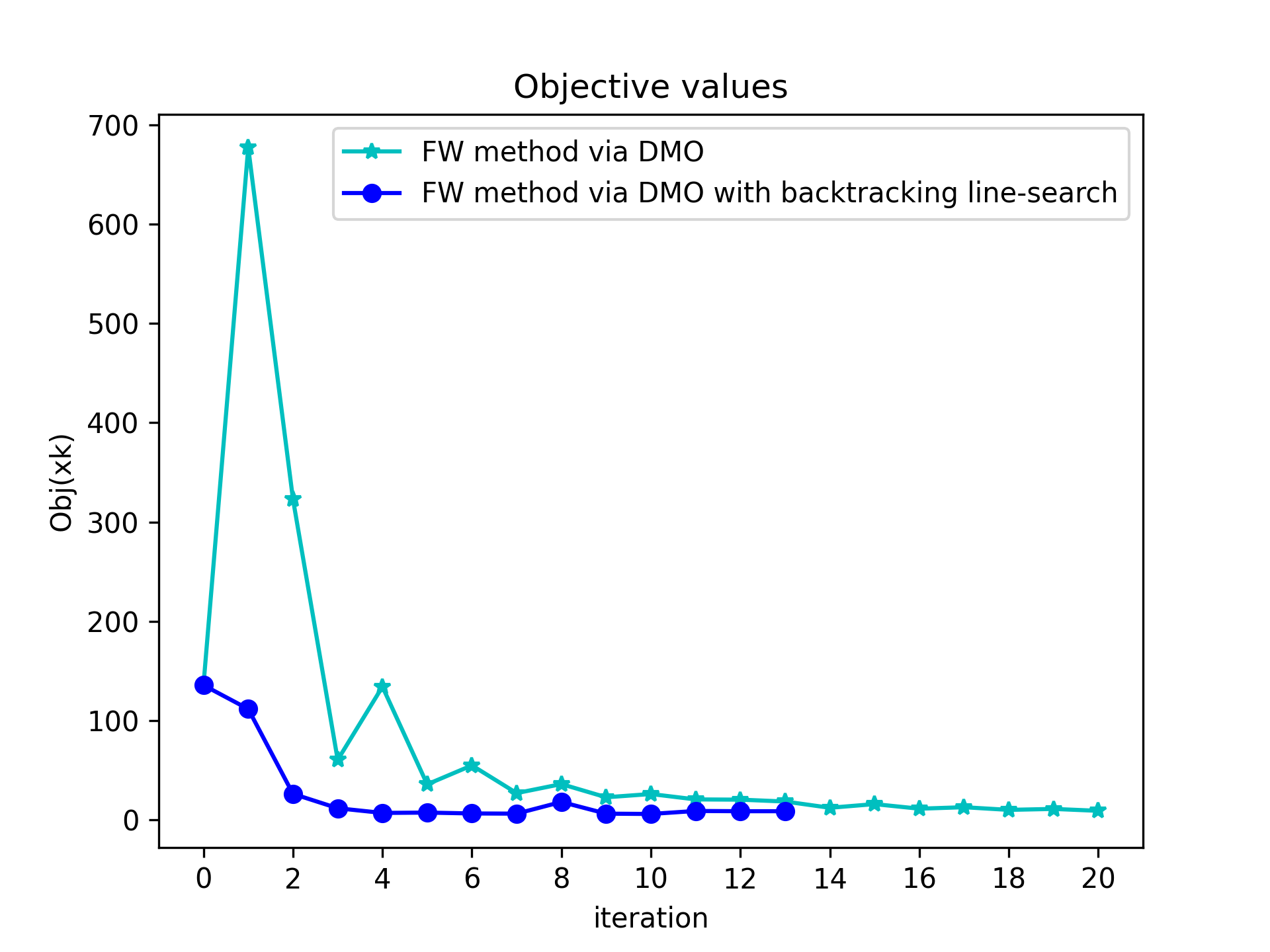}
\caption{Objective function values for the Frank-Wolfe methods via DMO, without and with backtracking line-search.}
\end{figure}

\begin{figure}[H]
\centering
\includegraphics[width=.6\textwidth]{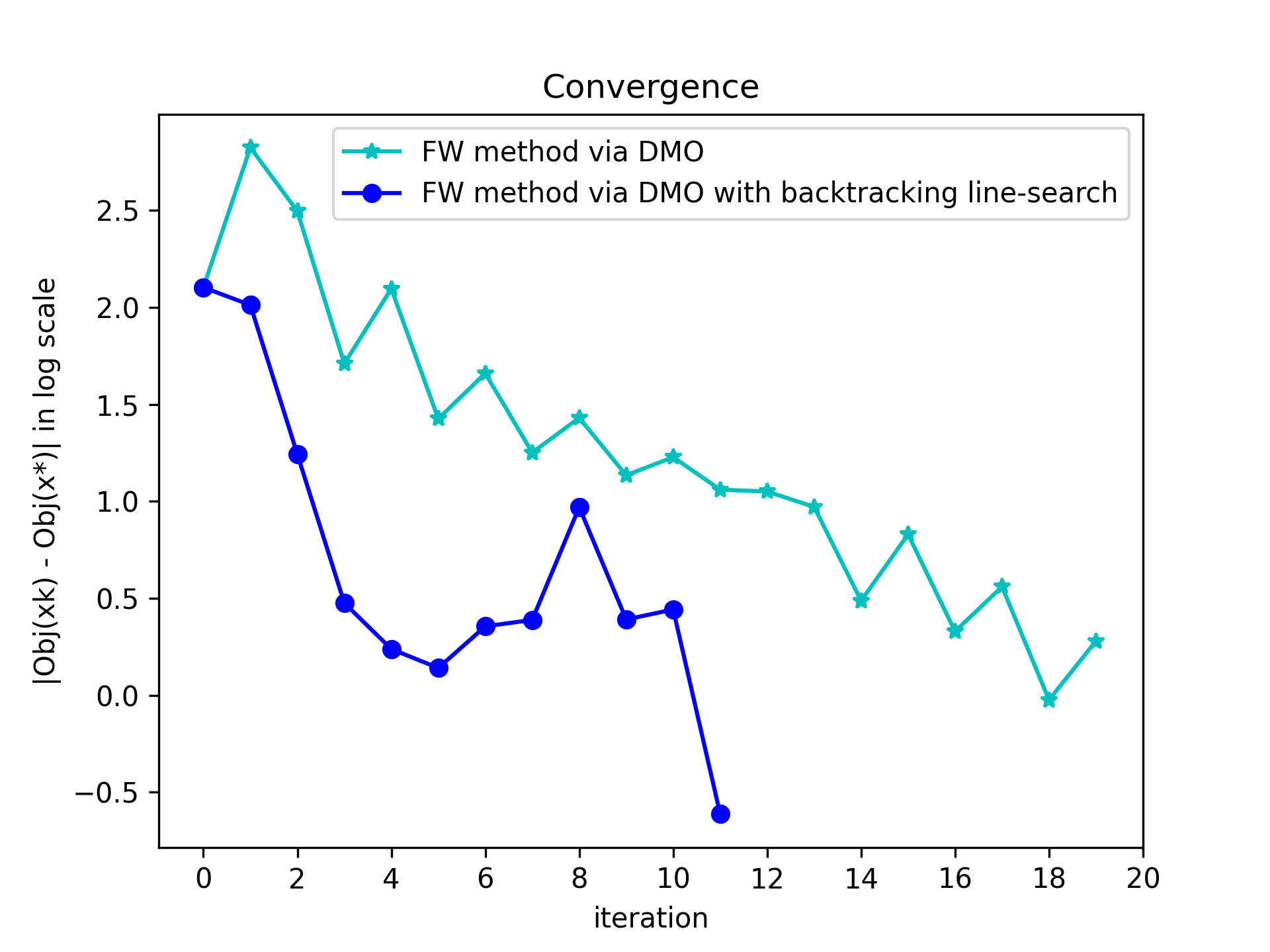}
\caption{Convergence plot for the Frank-Wolfe methods via DMO, without and with backtracking line-search.}
\end{figure}

In Fig. 4 and Fig. 5, it is apparent that the FW method via DMO with backtracking line search greatly reduces the objective function value at every iteration. The iteration times needed for convergence are also greatly reduced. This can be partially explained by the nature of the method because it ensures that the objective function value decreases sufficiently with the chosen step size, rather than a possible increase because of a larger step size. If the condition is not satisfied, the step size is reduced and the process is repeated until the condition is met. By using a backtracking line-search method to select the step size in the Frank-Wolfe algorithm, it can make better progress toward the solution with fewer iterations, leading to faster convergence.

To verify Top-g+ optimal visiting method we propose converges faster, we compared the result of Alg. 3 and Alg. 4 (in Appendix).  Blue line denotes the result of FW method via Top-g+ visiting while purple line denotes the result of FW method via Top-g+ optimal visiting.

\begin{figure}[H]
\centering
\includegraphics[width=.6\textwidth]{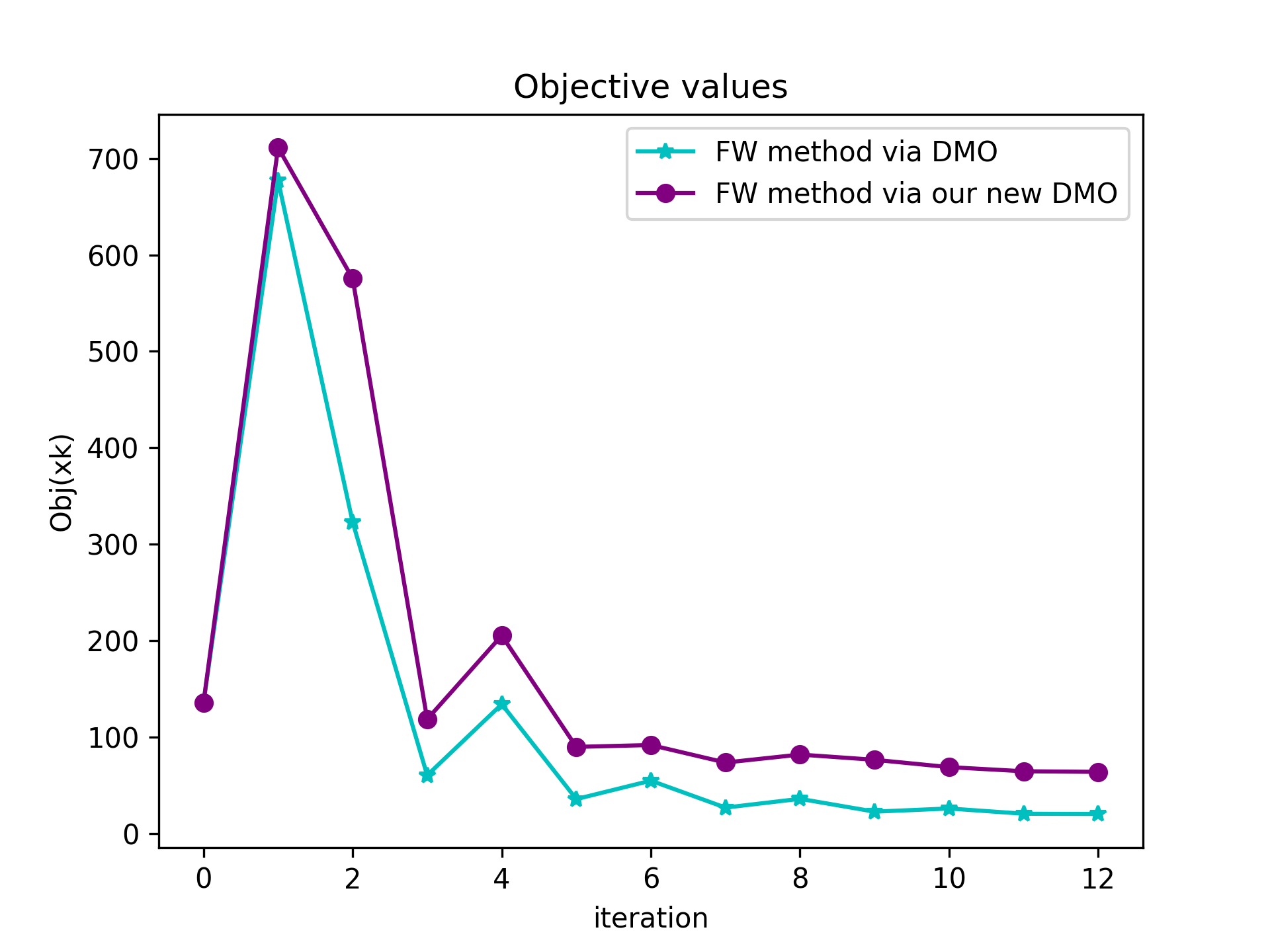}
\caption{Objective function values for the Frank-Wolfe methods via original DMO and our modified DMO.}
\end{figure}

\begin{figure}[H]
\centering
\includegraphics[width=.6\textwidth]{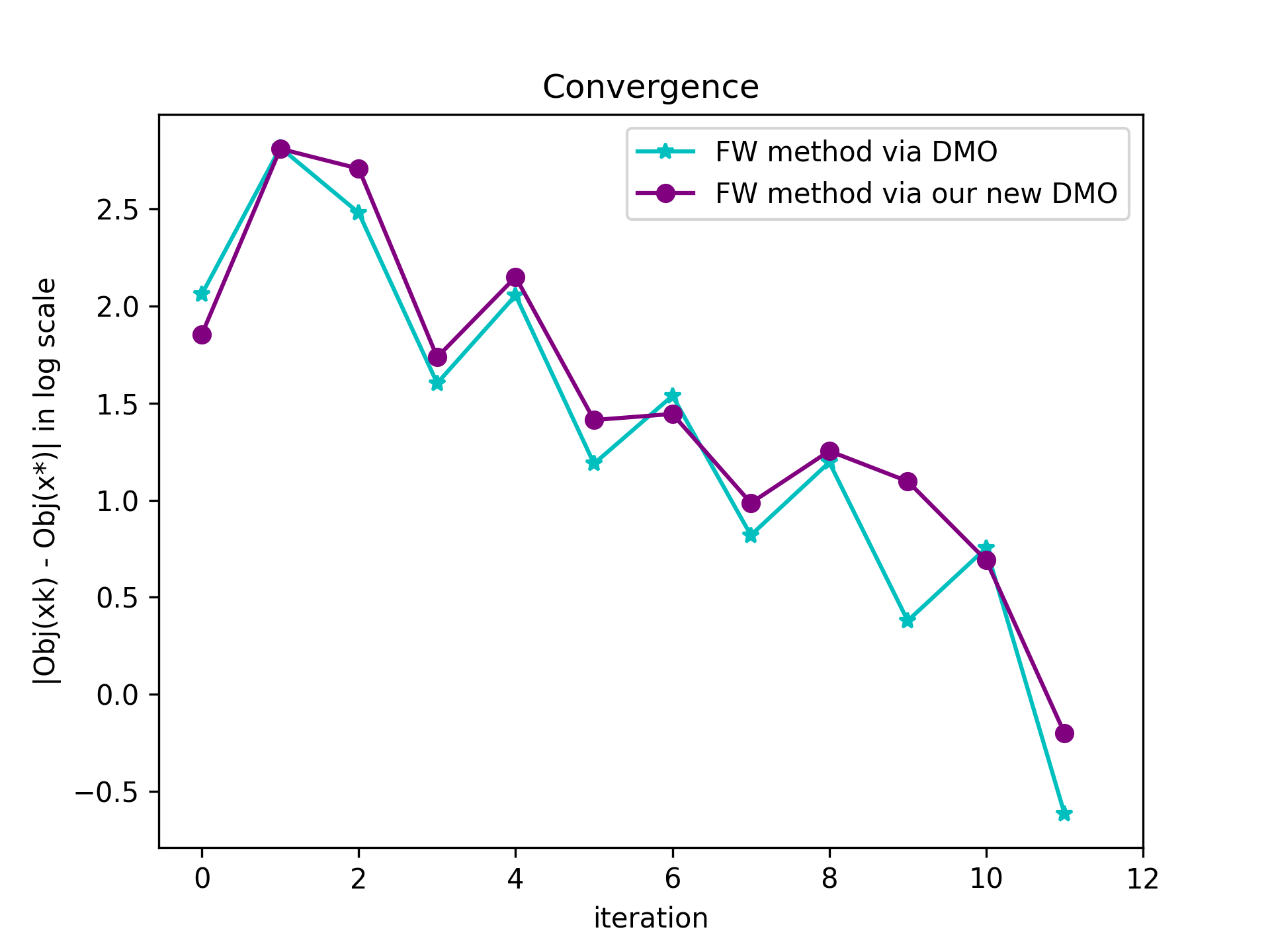}
\caption{Convergence plot for the Frank-Wolfe methods via original DMO and our modified DMO.}
\end{figure}

The comparison of objective values and convergence is shown in Fig. 6 and Fig. 7. Compared with FW method via the original DMO method (Top-g+ visiting), our new DMO method (Top-g+ optimal visiting) converges slightly faster in some points while the convergence rate can even be slower in other points. One possible reason is that the expected quantity of support sets we want to compare is not enough. Without adequate support sets to compare, we can not guarantee to find the optimal support set. However, the computational complexity will dramatically increase if we expand the search area for support sets.

\section{Conclusion}
In this report, we first introduced a kind of graph-structured convex optimization (GSCO) problem, which is relatively not easy to solve using traditional methods because the solution should satisfy the graph structure. It will always be a huge pain to obtain the optimal support set by traversing every group using LMO. Then we introduced a new DMO method to get an approximate IPO, which has more feasibility in some cases. In our experiments, we compared the objective function values for the FW method via DMO, the random PGD method, and the best PGD by comparing all possible cases. And we got a relatively good result: the outputs via the DMO method are very close to the Best PGD method. During the experiments, we also found that the authors of this paper didn't use a proper graph-like dataset, so they got an extremely perfect result. In the end, to deal with the shortcomings in the original algorithm, we proposed two improvements. The backtracking line-search method effectively reduced the number of iterations, while our new DMO method (Top-g+ optimal visiting) did not make many contributions. We notice that the improvement of convergence rate is not significant, but the running time will dramatically increase if we want Top-g+ optimal visiting to have better performance. Our experiments indicate that it is not reasonable to improve convergence rate at the huge cost of running time. However, at several points Top-g+ optimal visiting method indeed improves the convergence rate. The result proves that the optimal support set exists, but the method to find it efficiently still needs more studies.

\section*{References}
{
\small

[1] Zhou, B., and Sun, Y., “Approximate Frank-Wolfe Algorithms over Graph-structured Support Sets.” International Conference on Machine Learning, pp. 27303-27337. PMLR, 2021.
}

[2] Demyanov, Vladimir and Rubinov, Aleksandr. "Approximate Methods in Optimization Problems". Elsevier, 1970. 

[3] Dunn, J. C. and Harshbarger, S. Conditional gradient algorithms with open loop step size rules. Journal of Mathematical Analysis and Applications, 62(2):432–444, 1978.

[4] Frank, M., Wolfe, P., et al. "An algorithm for quadratic programming". Naval research logistics quarterly, 3(1-2): 95–110, 1956.

[5] Goldstein, Allen A. “On steepest descent”. Journal of the Society for Industrial and Applied Mathematics, Series A: Control, 1965.

[6] Jaggi, M. Revisiting Frank-Wolfe: Projection-free sparse convex optimization. In Proceedings of the 30th international conference on machine learning, pp. 427–435, 2013.

[7] Lacoste-Julien, Simon. "Convergence rate of Frank-Wolfe for non-convex objectives". arXiv preprint arXiv:1607.00345, 2016.

[8] Locatello, F., Yurtsevert, A., Fercoq, O., and Cevhert, V. "Stochastic Frank-Wolfe for composite convex minimization". In Proceedings of the 31st International Conference on Neural Information Processing Systems, volume 32, 2019.

[9] Pedregosa, Fabian and Askari, Armin and Negiar, Geoffrey and Jaggi, Martin (2018) "Step-Size Adaptivity in Projection-Free Optimization". arXiv:1806.05123

[10] Pedregosa, F., Negiar, G., Askari, A., and Jaggi, M. "Linearly convergent Frank-Wolfe with backtracking line-search". In International Conference on Artificial Intelligence and Statistics, pp. 1–10. PMLR, 2020.

[11] Armijo, L. "Minimization of functions having Lipschitz continuous first partial derivatives". Pacific Journal of Mathematics, 16(1), 1-3, 1966.

\clearpage
\section{Appendix}

\begin{algorithm}[H]
\caption{\textsc{Original DMO with $\delta=\sqrt{1/\lceil s/g\rceil}$} approximation guarantee}
\begin{algorithmic}[1]
 \STATE \textbf{Input}: underlying graph $\mathbb{G}(\mathbb{V},\mathbb{E})$, sparsity $k$, number of CCs $g$, input vector $\bm z$
 \STATE Sort entries of $\bm z$ by magnitudes such that $|z_{\tau_1}| \geq |z_{\tau_2}| \geq \ldots \geq |z_{\tau_g}|\geq |z_{\tau_{g+1}}|$ 
\STATE $I_g = [\tau_1,\tau_2,\ldots,\tau_g], S = I_g$
\STATE $\bm c = \bm 0$ \algorithmiccomment{Initially, all nodes have same connected component ID}
\STATE $i=1$ // Tracking the ID of connected component
\FOR{$v \in S$}
\STATE $c_{v} = i$ // Node $v$ has a component ID $i$
\STATE $i = i +1$ 
\ENDFOR
\STATE $\mathbb{F}=\emptyset$ // Keep edges that are in $g$ components
\IF{$|S| = s$} 
\STATE \textbf{Return} $S$ // We assume $g \leq s$
\ENDIF
\FOR{$(u,v) \in \mathbb{E}$}
\IF{$c_{u} == 0$ and $c_{v} \ne 0$}
\STATE $S = S \cup \{u\}$
\STATE $\mathbb{F} = \mathbb{F} \cup (u,v)$
\STATE $c_u = c_v$ // $u$ is added to $c_v$-th component
\ENDIF
\IF{$|S| = s$}
\STATE \textbf{Return} $S$
\ENDIF
\IF{$c_{u} \ne 0$ and $c_{v} == 0$}
\STATE $S = S \cup \{v\}$
\STATE $\mathbb{F} = \mathbb{F} \cup (u,v)$
\STATE $c_v = c_u$ // $v$ is added to $c_u$-th component
\ENDIF
\IF{$|S| = s$}
\STATE \textbf{Return} $S$
\ENDIF
\ENDFOR
\end{algorithmic}
\label{algo:original heuristic-dmo}
\end{algorithm}

\begin{algorithm}[h]
\caption{\textsc{new DMO with $\delta=\sqrt{1/\lceil s/g\rceil}$} approximation guarantee}
\begin{algorithmic}[1]
 \STATE \textbf{Input}: underlying graph $\mathbb{G}(\mathbb{V},\mathbb{E})$, number of CCs $g$, input vector $\bm z$, expected quantity of support sets $\theta$ 
 \STATE Sort entries of $\bm z$ by magnitudes such that $|z_{\tau_1}| \geq |z_{\tau_2}| \geq \ldots \geq |z_{\tau_g}|\geq |z_{\tau_{g+1}}|$ 
\STATE $I_g = [\tau_1,\tau_2,\ldots,\tau_g], S = I_g$
\STATE $\bm c = \bm 0$ \algorithmiccomment{Initially, all nodes have same connected component ID}
\STATE $i=1$ // Tracking the ID of connected component
\FOR{$v \in S$}
\STATE $c_{v} = i$ // Node $v$ has a component ID $i$
\STATE $i = i +1$ 
\ENDFOR
\STATE $\mathbb{F}=\emptyset$ // Keep edges that are in $g$ components
\IF{$|S| = s$} 
\STATE \textbf{Return} $S$ // We assume $g \leq s$
\ENDIF

\STATE $j =0$
\STATE $\mathbb{E}_0=\mathbb{E}$
\STATE $S_0=S$
\STATE $n=(s-g)*\theta$
\FOR{$t = 0,1,\ldots,n$}
\STATE randomly select $ (u,v) \in \mathbb{E}$
\IF{$c_{u} == 0$ and $c_{v} \ne 0$}
\STATE $S = S \cup \{u\}$
\STATE $\mathbb{F} = \mathbb{F} \cup (u,v)$
\STATE $c_u = c_v$ // $u$ is added to $c_v$-th component
\ENDIF
\IF{$|S| = s$}
\STATE  $out[j]=\left \|z _S\right\|_2$
\STATE $S_{list}[j]=S$ // record possible support sets
\STATE $j =j+1$
\STATE $\mathbb{E}=\mathbb{E}_0$
\STATE $S=S_0$
\ENDIF
\IF{$c_{u} \ne 0$ and $c_{v} == 0$}
\STATE $S = S \cup \{v\}$
\STATE $\mathbb{F} = \mathbb{F} \cup (u,v)$
\STATE $c_v = c_u$ // $v$ is added to $c_u$-th component
\ENDIF
\IF{$|S| = s$}
\STATE  $out[j]=\left \|z _S\right\|_2$ 
\STATE $S_{list}[j]=S$ // record possible support sets
\STATE $j =j+1$
\STATE $\mathbb{E}=\mathbb{E}_0$
\STATE $S=S_0$
\ENDIF
\ENDFOR
\STATE d=index(max(out))
\STATE $S=S_{list}[d]$
\STATE \textbf{Return} $S$
\end{algorithmic}
\label{algo:heuristic-dmo}
\end{algorithm}

\end{document}